\documentclass[12pt,draft,leqno]{article}
\usepackage{amssymb, eucal, latexsym}

\newtheorem{theorem}{Theorem}[section]
\newtheorem{lm}[theorem]{Lemma}
\newtheorem{exa}[theorem]{Example}

\newtheorem{cor}[theorem]{Corollary}
\newtheorem{pro}[theorem]{Proposition}
\newtheorem{defi}[theorem]{Definition}

\newtheorem{rem}[theorem]{Remark}

\newtheorem{question}[theorem]{Question}

\def\p{\varphi}
\def\a{\alpha}
\def\b{\beta}

\def\ep{\varepsilon}
\def\g{\gamma}
\def\GA{\Gamma}

\def\l{\lambda}

\def\s{\sigma}
\def\SI{\Sigma}

\def\OM{\Omega}

\def\lra{\longrightarrow}

\def\sbe{\subseteq}
\def\spe{\supseteq}
\def\stm{\setminus}
\def\ems{\emptyset}
\def\nes{\neq\emptyset}

\def\ex{\exists}
\def\fa{\forall}

\def\ap{^{\,\prime}}
\def\inv{^{-1}}
\def\st{\ |\ }

\def\nin{\not\in}

\def\AA{{\cal A}}
\def\BB{{\cal B}}

\def\EE{{\cal E}}

\def\NN{{\cal N}}

\def\PP{{\cal P}}
\def\QQ{{\cal Q}}

\def\TT{{\cal T}}

\def\UU{{\cal U}}
\def\VV{{\cal V}}
\def\WW{{\cal W}}

\def\ZHC{{\bf Stone}}
\def\Bool{{\bf Bool}}

\def\2{\mbox{{\bf 2}}}
\def\3{\mbox{{\bf 3}}}

\def\int{\mbox{{\rm int}}}

\def\cl{\mbox{{\rm cl}}}

\def\doc{\hspace{-1cm}{\em Proof.}~~}
\def\sq{\hspace*{\fill} \hbox{\vrule\vbox{\hrule\phantom{o}\hrule}\vrule}}
\def\sqs{\sq \vspace{2mm}}

\def\bk{\bar{k}}
\def\bn{\bar{n}}

\def\RRRR{\mathbb{R}}
\def\NNNN{\mathbb{N}}
\def\DDDD{\mathbb{D}}
\def\QQQQ{\mathbb{Q}}

\def\IIII{\mathbb{I}}
\def\PPPP{\mathbb{P}}

\title{{\LARGE\bf
An internal topological characterization of the subspaces of Eberlein compacta and related compacta -- II}\\
\vspace{0.35cm}
{\large\bf Georgi D. Dimov}\thanks{This paper was supported by the
project no. 140/2014 $``$Function Spaces and Dualities" of the
Sofia
University $``$St. Kl. Ohridski".
}\\
\vspace{0.25cm}
 {\footnotesize Dept. of Math. and
Informatics, Sofia University,  5 J. Bourchier Blvd., 1164 Sofia,
Bulgaria}
}

\author{}

\date{}

\begin{document}
\maketitle
\begin{abstract}
We generalize a theorem of E. Michael and M. E. Rudin \cite[Theorem 1.1]{MR2} and a theorem of D. Preiss and P. Simon \cite[Corollary 7]{PrS}; we give, as well,
some partial answers to a recent question of A. V. Arhangel'ski\v{\i}.
\end{abstract}



\footnotetext[3]{{\footnotesize {\em E-mail address:}
gdimov@fmi.uni-sofia.bg}}


\baselineskip = \normalbaselineskip

\section{Introduction}


\bigskip

In the first part \cite{D4} of this paper, we supplied with proofs a part of the results announced (without any proofs) in our paper \cite{D1}. Here we will give the proofs of the remaining part of the assertions from \cite{D1} and we will obtain some new results as well.  The paper \cite{D1} was inspired by the following
 question of A. V. Arhangel'ski\v{\i} (posed in 1982): {\em $``$When does a space $X$ have a compactification $cX$ which is an Eberlein compact?"}, and by some related problems.
The  above question arose naturally from the following  theorem of
 A. V. Arhangel'ski\v{\i} proved in 1977 (\cite[Theorem 14]{A}):  {\em every metrizable space has a compactification which is an Eberlein compact}.
 Since the closed subspaces of Eberlein compacta are Eberlein compacta,
 the Arhangel'ski\v{\i} question  is equivalent to the following problem: {\em find an internal characterization of  the subspaces of Eberlein compacta}. Let's note that in 1974 H. P. Rosenthal \cite{Ro} gave an internal characterization of Eberlein compacta (see Theorem \ref{Rosenthal} below).  In \cite{D1,D4} we characterized internally the subspaces of Eberlein compacta, as well as the subspaces of Corson compacta, of uniform Eberlein compacta, of $n$-uniform Eberlein compacta (in the sense of B. A. Pasynkov \cite{P}) and of strong Eberlein compacta.

 This paper is organized as follows. {\em Section 2}\/ contains some preliminary results and some definitions which are indispensable for our exposition.
 In {\em Section 3}\/ we generalize (by our Theorem \ref{GMR}) a  theorem of E. Michael and M. E. Rudin \cite[Theorem 1.1]{MR2} (see Theorem \ref{MR2Th11} below). The proof of  Theorem \ref{GMR} is based on the  Michael-Rudin characterization of Eberlein compacta \cite[Theorem 1.4, (a)$\leftrightarrow$(e)]{MR1}
 and on our Theorem \ref{th87} in which, modifying Michael-Rudin's proof of \cite[Theorem 1.4, (a)$\leftrightarrow$(e)]{MR1}, we obtain a new characterization of uniform Eberlein compacta which is similar to that of Eberlein compacta given by E. Michael and M. E. Rudin.
  In {\em Section 4}\/ we generalize (by our Theorem \ref{GPS0}) a theorem of D. Preiss and P. Simon \cite[Corollary 7]{PrS} (see Theorem \ref{e6} below). For doing this, we introduce the notion of {\em $F_\tau$-embedded subspace}, where $\tau$ is a cardinal number. We give some sufficient conditions for $F_\tau$-embeddability (see Theorem \ref{Ftauemb}). The last theorem implies, in particular, that: (a) every metrizable dense subspace $Y$ of a space $X$ is $F_{\aleph_0}$-embedded in $X$, (b) if $X$ is a $T_4$-space, $Y\sbe X$ and $nw(Y)\le \tau$, then $Y$ is $F_\tau$-embedded in $X$. Using the ideas of the proof of Theorem \ref{GPS0}, we obtain some  conditions such that if a rim-compact space $X$ satisfies them, then its Freudenthal compactification $FX$ is an Eberlein compact or a strong  Eberlein compact (see Theorem \ref{fr} and its corollaries). Almost all of the results from Sections 3 and 4 were announced without any proofs in our paper \cite{D1}, and the remaining results from these sections are from \cite{D0}. {\em Section 5}\/  contains only new results. With them we give some partial answers to the following recent question of A. V. Arhangel'ski\v{\i}: {\em $``$Is it true that if a Tychonoff space $X$ has a $\s$-disjoint base then it can be embedded in an Eberlein compact?"}. We prove, in particular, that any $T_0$-space which has a $\s$-disjoint base consisting of clopen sets can be embedded in a zero-dimensional uniform Eberlein compact (see the more general Theorem \ref{th83}). We show, as well, that if the answer to the just mentioned question of A. V. Arhangel'ski\v{\i} is $``$Yes", then every Tychonoff space with a $\s$-disjoint base has a $\s$-point-finite base consisting of cozero-sets (see Corollary \ref{c899}).

We now fix the notation.

If $A$ is a set, we denote by $|A|$ its cardinality and by $exp(A)$ its power set.
If $\UU$ is a family of subsets of $A$ and $a\in A$, then we set $$\UU(a)=\{U\in\UU\st a\in U\}.$$
 If $(X,\tau)$ is a topological space and $M$ is a subset of $X$, we
denote by $\cl_{(X,\tau)}(M)$ (or simply by $\cl(M)$ or
$\cl_X(M)$) the closure of $M$ in $(X,\tau)$ and by
$\int_{(X,\tau)}(M)$ (or briefly by $\int(M)$ or $\int_X(M)$) the
interior of $M$ in $(X,\tau)$. The
set of all positive natural numbers is denoted  by  $\NNNN$,
  the real line (with its natural topology) -- by $\mathbb{R}$, and the subspaces $[0,1]$ and $\{0,1\}$ of $\mathbb{R}$ -- by $\IIII$ and $\DDDD$, respectively.
  As usual, $\omega=\NNNN\cup\{0\}$.
  The Stone-\v{C}ech compactification of a Tychonoff space $X$ is denoted by $\b X$.

Let $X$ be a dense subspace of a space $Y$ and $U\sbe X$. The {\em extension of $U$ in $Y$}, denoted by $Ex_YU$, is the set $Y\stm\cl_Y(X\stm U)$; recall that if $U$ is open in $X$ then $Ex_YU$ is the largest open subset of $Y$ whose trace on $X$ is $U$.

By a {\em neighborhood} $U$ of a point $x$ of a topological space $X$ we mean a neighborhood in the sense of Bourbaki (i.e., $x\in \int_X(U)$).

If $X$ is a topological space then we denote by $R(X)$ the set of all points of $X$ without compact neighborhoods.

If  $X$ is a topological space and $f:X\lra \IIII$ is a continuous function, then we write, as usual, $Z(f)=f\inv(0)$ ({\em the zero-set of $f$}\/) and $coz(f)=X\stm Z(f)$ ({\em the cozero-set of $f$}\/). We set $Z(X)=\{Z(f)\st
f:X\lra \IIII$ is a continuous function$\}$ and $Coz(X)=\{coz(f)\st
f:X\lra \IIII$ is a continuous function$\}$.

We  denote by $\QQQQ$ the set of all  rational numbers, by $\PPPP$ the set of all  irrational numbers, and by $\QQQQ_1$ the set of all positive rational numbers less than 1.

 We  denote by  $CO(X)$
the collection  of all  clopen (= closed and open) subsets of $X$.

We denote by
$S:\Bool\lra \ZHC$ the
Stone duality functor between the category $\Bool$ of Boolean algebras and Boolean
homomorphisms and the category $\ZHC$ of compact
zero-dimensional Hausdorff spaces  and
continuous maps (see, e.g., \cite{kop89}).

All notions and notation which are not explained here can be found in \cite{E,kop89}.

{\bf By a $``$space" we  will mean a $``$topological $T_0$-space".}

\section{Preliminaries}
%

Let  $\GA$ be an index set and let $\IIII^\GA$ (respectively, $\DDDD^\GA$) be the Cartesian product of $|\GA|$
copies of $\IIII$ (respectively, $\DDDD$) with the Tychonoff topology. We set
$$\SI_*(\IIII,\GA)=\{ x\in \IIII^\GA\st (\fa \ep> 0) (|\{\g\in\GA\st x(\g)\ge\ep\}|<\aleph_0) \},$$
$$\s(\DDDD,\GA) =\{ x\in \DDDD^\GA\st |\{\g\in\GA\st x(\g)=1\}|<\aleph_0\},$$
and the topology on these subsets of $\IIII^\GA$ is the subspace topology. Obviously, $\s(\DDDD,\GA)\sbe\SI_*(\IIII,\GA)$.

Let $\UU$ be a family of subsets of a set $X$. Recall that  by an {\em order}\/  of the family  $\UU$  (denoted by $ord(\UU)$) we mean the largest integer $n$ such that the family $\UU$ contains $n+1$ sets with non-empty intersection, or the $``$infinite number" $\infty$, if no such integer exists. Thus, if $ord(\UU)=n\in\omega$, then any $n+2$
 members of
the family $\UU$ have an empty intersection.
 The family $\UU$ is called {\em boundedly point-finite}\/  (\cite{FMS}) if
 it has a finite order. If  the family $\UU$ is a union of countably many families $\UU_n$, where $n\in\NNNN$, we will say that $\UU$ is a {\em $\s$-point-finite}\/ (respectively, {\em $\s$-boundedly point-finite}\/) {\em  family}\/  if all families $\UU_n$ are point-finite (respectively, boundedly point-finite).
 The family $\UU$ is said to be {\em $T_0$-separating}\/ if, whenever $x\neq y$ are in $X$, then there exists $U\in\UU$ such that $|U\cap\{x,y\}|=1$; in this case we will also say  that the family $\UU$ {\em $T_0$-separates $X$}.
Finally, when $X$ is a topological space and $Y\sbe X$, we say that
a family $\UU$ of subsets of $X$  {\em F-separates} $Y$
(see \cite{MR1}) if, whenever $x\neq y$ are in $Y$,
then there is some $U\in\UU$ such that $x\in U$ and $y\nin \cl_X(U)$, or vice versa; when $Y=X$, we will say that $\UU$ is an {\em F-separating family in} $X$.

 A compact Hausdorff space is called an {\em Eberlein
compact}\/ (briefly, EC), if it is homeomorphic to
a weakly compact  (i.e., compact in the weak topology) subset of a Banach space (\cite{L}).
D. Amir and J. Lindenstrauss proved that a compact
space is an Eberlein compact if and only if it can be embedded in
$\SI_*(\IIII,\GA)$ for some index set $\GA$  (see \cite[Theorem 1]{AL}). An internal characterization of Eberlein compacta was given by H. P. Rosenhthal \cite{Ro}.

\begin{theorem}\label{Rosenthal}{\rm (\cite{Ro})}
A compact Hausdorff space $X$ is an
Eberlein compact if and only if it has a $\s$-point-finite $T_0$-separating
collection of cozero-sets.
 \end{theorem}

A compact Hausdorff space that is homeomorphic to
a weakly compact subset of a Hilbert space is called a {\em uniform Eberlein compact}\/ (briefly, UEC); this class of spaces was introduced by Benyamini and Starbird in \cite{BS}. Let us recall the following result:

\begin{theorem}\label{argf}{\rm (\cite{BRW})}
Let $X$ be a compact Hausdorff space. Then  $X$ is a uniform Eberlein compact iff
$X$ has a $\s$-boundedly point-finite, $T_0$-separating family of cozero sets.
 \end{theorem}

A compact space $X$ is called a {\em strong Eberlein compact}\/ (briefly, SEC) (\cite{Sim,BRW,Ro}) if it can be embedded in $\s(\DDDD,\GA)$ for some set $\GA$.
Equivalently, a compact space is a SEC iff it has a point-finite, $T_0$-separating family of clopen  sets.
K. Alster \cite{Alster} proved that a space is a SEC iff it is a scattered EC.

\begin{defi}\label{Urepr}{\rm (\cite{D1,D2})}
\rm
 Let $X$ be a space and $V$ a subset of it. If there
exists a collection $U (V) =\{U_n(V)\st n\in\NNNN\}$ such that $V =\bigcup\{U_n(V) \st
n\in\NNNN\}$ and $U_n(V)\sbe U_{n+1}(V)$, $U_{2n-1}(V)\in Z(X)$, $U_{2n}(V)\in Coz(X)$
for every $n\in\NNNN$, then we  say that the set $V$ is {\em U-representable}\/
and the collection $U(V)$ is a {\em U-representation of} $V$. If $\a$ is
a family of subsets of $X$ and, for every $V\in\a$, $U(V)$ is a
U-representation of $V$, then the family $U(\a)= \{U(V) \st V \in \a\}$
is called a {\em U-representation of} $\a$.
\end{defi}

We will need the following lemma:

\begin{lm}\label{Uf}{\rm (\cite{D2})}
 Let $X$ be a space and $U(V)=\{U_n(V)\st n\in\NNNN\}$ be a U-representa\-ti\-on of a subset $V$ of $X$. Then there exists a continuous function $f:X\lra\IIII$ such that
 $V=coz(f)$ and $f\inv([\frac{1}{2^{n-1}},1])=U_{2n-1}(V)$, for every $n\in\NNNN$.
\end{lm}

The definition given below was inspired by the results
of B. Efimov and G. \v{C}ertanov \cite{EC}:

\begin{defi}\label{alsubb}{\rm (\cite{D1})}
\rm
 A  family $\a$ of subsets of a space $X$ is said to be an {\em almost subbase
of} $X$ if there exists a U-representation $U(\a)$ of $\a$
 such that the family $\a\cup \{X\stm U_{2n-1}(V)\st V\in\a, n\in\NNNN\}$ is a subbase of $X$.
\end{defi}

\begin{rem}\label{rem00}
\rm
Every almost subbase $\a$ of a space $X$ consists of cozero-sets because a subset of $X$ is U-representable iff
it is a cozero-set. Obviously, a U-representable set $V$ can have many different
U-representations.

It is easy to see that a space has an almost subbase iff it is completely regular (see \cite{D1,D2}).
\end{rem}

\begin{defi}\label{strpf}{\rm (\cite{Pe1,EP,D1})}
\rm
A family $\a$ of subsets of a set $X$ is called {\em strongly point-finite}\/  if every subfamily $\mu$ of $\a$ with $|\mu|=\aleph_0$ contains
a finite subfamily $\mu\ap$ with empty intersection.
A family of subsets of a set $X$ is called  {\em $\s$-strongly point-finite}\/ if it is a union of countably many  strongly point-finite families.
\end{defi}

 Obviously,  every strongly point-finite  family is a point-finite  family; the converse is not true.
 
 Recall that a family $\{A_s\st s\in S\}$ of subsets of a set $X$ is {\em star-finite}\/ if for every $s_0\in S$ the set $\{s\in S\st A_s\cap A_{s_0}\nes\}$ is finite (see, e.g., \cite{E}). Clearly, {\em any star-finite family is strongly point-finite.} The converse is not true. For example, in $\mathbb{R}^2$, the family $\AA=\{\{\frac{1}{n}\}\times [0,1]\st n\in\NNNN\}\cup \{[0,1]\times\{0\}\}$ is strongly point-finite but it is not star-finite.

\begin{rem}\label{rem0}
\rm
In \cite{D4} we noted that
if $\a$ is an almost subbase of a $T_0$-space then the family $$\a\ap=\a\cup\{f_V\inv((r,1])\st V\in\a, r\in\QQQQ_1\}$$ (where, for any $V\in\a$, $f_V$ is the function constructed in Lemma \ref{Uf} on the base of the given U-representation $U(V)$ of $V$) is  an F-separating family and an almost subbase. Moreover, when $\a$ is a $\s$-point-finite (respectively, $\s$-strongly point-finite; $\s$-boundedly point-finite) family then $\a\ap$ has the same property.
\end{rem}

The next lemma will be used frequently:

\begin{lm}\label{Ex}{\rm (\cite{D0,D4})}
 Let $X$ be a dense subspace of a space $Y$ and $\g$ be a strongly point-finite (respectively, boundedly point-finite) family of open subsets of $X$. Then the family $$Ex_Y\g=\{Ex_YU\st U\in\g\}$$ is a strongly point-finite (respectively, boundedly point-finite) family of open subsets of $Y$.
\end{lm}

\begin{theorem}\label{unifthE}{\rm (\cite{D4})}
Let $X$ be a compact Hausdorff space.  Then $X$ is an EC (respectively,  UEC) iff
 $X$ has a $\s$-strongly point-finite  almost subbase  (respectively,  $\s$-boundedly point-finite   almost subbase);
  $X$ is a strong Eberlein compact iff it has a strongly point-finite family  $\a$ consisting of clopen sets such that $\a\cup\{X\stm V\st V\in\a\}$ is a subbase of $X$.
 \end{theorem}

\begin{defi}\label{espa}{\rm (\cite{D1})}
\rm
  A space $X$ is called an {\em Eberlein space}\/ (brieily, {\em E-space}\/) (respectively,
{\em UE-space}\/) if $X$ has a $\s$-strongly point-finite (respectively, $\s$-boundedly
point-finite) almost subbase. A space $X$ is said to be a {\em SE-space}\/ if it has a strongly
point-finite family $\a$ of clopen subsets such that $\a\cup \{X\stm V\st  V\in\a\}$ is a subbase of $X$.
\end{defi}

\begin{theorem}\label{MainE}{\rm (\cite{D1,D4})}
  A space X has a compactification $cX$ which is an Eberlein compact (respectively,  uniform Eberlein compact; strong Eberlein compact) iff it is an E-space
 (respectively, UE-space; SE-space).
\end{theorem}

\begin{cor}\label{subsp}{\rm (\cite{D1,D4})}
 A space $X$ is a subspace of an Eberlein compact (respectively,  UEC; SEC)  iff it is an E-space (respectively,  UE-space; SE-space).
 \end{cor}



\section{On a theorem of E. Michael and M. E. Rudin}

In \cite{MR2}, E. Michael and M. E. Rudin proved the following theorem:

\begin{theorem}\label{MR2Th11}{\rm (\cite[Theorem 1.1]{MR2})}
If $X$ is compact Hausdorff, and if $X=X_1 \cup X_2$ with $X_1$ and $X_2$ metrizable, then $\cl(X_1)\cap\cl(X_2)$ is metrizable and
$X$ is an EC.
\end{theorem}

We are now going to generalize it. We will need some preliminary results.

Let's first recall the following lemmas from \cite{MR1}:

\begin{lm}\label{MR1Lm32}{\rm (\cite[Lemma 3.2]{MR1})}
Let $(X,\TT)$ be Tychonoff, $\l=|X|$ and $U\in\TT$. Then there is an open cover $\{V_\a\st \a<\l\}$ of $U$, and a closed (in $X$) cover $\{S_{\a,m}\st\a<\l,m\in\NNNN \}$ of $U$ , such that:

\smallskip

\noindent(a) $S_{\a,m}\sbe V_\a\sbe U$ for all $\a,m$,

\smallskip

\noindent(b) If $\a>\b$, then $S_{\a,m}\cap V_\b=\ems$ for all $m$.
\end{lm}

\begin{lm}\label{MR1Lm33}{\rm (\cite[Lemma 3.3]{MR1})}
Let $X$ be a space, $\UU$ be an F-separating  open cover of $X$, and $A,B$ be disjoint compact subsets of $X$. Then there is a finite $\EE\sbe \UU$ which F-separates any $x\in A$ from any $x\ap\in B$.
\end{lm}

In \cite{MR1}, E. Michael and M. E. Rudin proved the following theorem (since we will use later on the proof of this theorem and the notation introduced there, we will recall  here (for convenience of the reader) the Michael-Rudin proof of this theorem):

\begin{theorem}\label{MR1Th14}{\rm (\cite[Theorem 1.4]{MR1})} A compact Hausdorff space $X$ is an EC iff it has an F-separating
$\sigma$-point-finite family of open subsets.
\end{theorem}

\doc ($\Rightarrow$) By the Rosenthal Theorem \ref{Rosenthal}, $X$ has a $\s$-point-finite $T_0$-separating
collection $\UU$ of cozero-sets. Fix, for every $U\in\UU$, a function $f_U:X\lra [0,1]$ such that $U=f\inv((0,1])$. Then the family $\UU\ap=\{f_U\inv((\frac{1}{n},1])\st n\in\NNNN , u\in\UU\}$ is the required F-separating
$\sigma$-point-finite family of open subsets.

\smallskip

\noindent($\Leftarrow$) Let $X$ have an F-separating family $\UU=\bigcup\{\UU_n\st n\in\NNNN \}$ of open subsets, where all families $\UU_n$, $n\in\NNNN $, are point-finite. We may suppose that $\UU_n\sbe \UU_{n+1}$ and $\UU_n$ covers $X$, for all $n\in\NNNN $.

Let $\l=|X|$ . For every $U\in\UU$ choose, using Lemma \ref{MR1Lm32}, an open cover  $\{V_\a(U)\st \a<\l\}$ of $U$, and a closed (in $X$) cover $\{S_{\a,m}(U)\st\a<\l,m\in\NNNN \}$ of $U$, such that:

\smallskip

\noindent(a${}_U$) $S_{\a,m}(U)\sbe V_\a(U)\sbe U$ for all $\a,m$,

\smallskip

\noindent(b${}_U$) if $\a>\b$, then $S_{\a,m}(U)\cap V_\b(U)=\ems$ for all $m$.

\smallskip

For every $U\in\UU$, for every $\a<\l$ and for every $m\in\NNNN $, we shall define a positive integer $k=k_{\a,m}(U)$ such that
\begin{equation}\label{42}
\mbox{if } x,x\ap\in X, \mbox{ if } \UU_k(x)=\UU_k(x\ap) \mbox{ and if } x\in S_{\a,m}(U), \mbox{ then } x\ap\in V_\a(U).
\end{equation}

To define $k$, let $A=S_{\a,m}(U)$ and $B=X\stm V_\a(U)$. Then, by Lemma \ref{MR1Lm33}, there exists a finite subfamily $\EE$ of $\UU$ which F-separates any $x\in A$ from any $x\ap\in B$. Pick $k$ so that $\EE\sbe\UU_k$; this $k$ satisfies (\ref{42}).

Let $U\in\UU$, $\a<\l$, $m\in\NNNN$ and let $k=k_{\a,m}(U)$. Define
$$H_{\a,m}(U)=\bigcup\{\bigcap\UU_k(x)\st x\in S_{\a,m}(U)\}.$$
Then $H_{\a,m}(U)$ is an open subset of $X$ and $H_{\a,m}(U)\spe S_{\a,m}(U)$. Choose an open $F_\s$-subset $W_{\a,m}(U)$ of $X$ such that
$$S_{\a,m}(U)\sbe W_{\a,m}(U)\sbe H_{\a,m}(U)\cap V_\a(U)\sbe U.$$

Define
$$\WW=\{W_{\a,m}(U)\st U\in\UU, \a<\l, m\in\NNNN\}.$$
We will show that $\WW$ is a $\s$-point-finite $T_0$-separating family of $F_\s$-subsets of $X$, thereby completing our proof (using the Rosenthal Theorem \ref{Rosenthal}).

\smallskip

\noindent(A.) {\em $\WW$  $T_0$-separates $X$:}
Let $x,y\in X$ and $x\neq y$; then $\ex U\in\UU$ such that $|U\cap\{x,y\}|=1$. Let, e.g., $x\in U$ and $y\nin U$. Then $\ex \a<\l$ and $\ex m\in\NNNN$ such that $x\in
S_{\a,m}(U)$. Hence $x\in W_{\a,m}(U)$ and $y\nin W_{\a,m}(U)$.

\smallskip

\noindent(B.) {\em $\WW$ is $\s$-point-finite:} For $m,n,k\in\NNNN$, let
$$\WW_{m,n,k}=\{W_{\a,m}(U)\st U\in\UU_n, \a<\l, k=k_{\a,m}(U)\}.$$
Then, obviously, $\WW=\bigcup\{\WW_{m,n,k}\st m,n,k\in\NNNN\}$, so that we only need to show that the family $\WW_{m,n,k}$ is point-finite for every $m,n,k\in\NNNN$.

Suppose some $\WW_{m,n,k}$ were not point-finite. Then there are distinct $W_{\a_j,m}(U_j)$, $j\in\NNNN$, in $\WW_{m,n,k}$ all containing some point $x^*\in X$. Then
$x^*\in H_{\a_j,m}(U_j)$ for all $j\in\NNNN$, so there are $x_j\in S_{\a_j,m}(U_j)$, $j\in\NNNN$, such that $x^*\in\bigcap\UU_k(x_j)$ for all $j\in\NNNN$. Hence $\UU_k(x^*)\spe\UU_k(x_j)$ for all $j\in\NNNN$. Since $\UU_k(x^*)$ is finite, we may assume, by passing to a subsequence, that $\UU_k(x_i)=\UU_k(x_j)$ for all $i,j\in\NNNN$. Now, by (\ref{42}), it follows that $x_i\in V_{\a_j}(U_j)$ for all $i,j\in\NNNN$. Then
\begin{equation}\label{43}
x_i\in S_{\a_i,m}(U_i)\cap V_{\a_j}(U_j) \mbox{ for all } i,j\in\NNNN.
\end{equation}
We now distinguish two cases to obtain our contradiction:

\smallskip

\noindent{\em Case 1.}  {\em There are $i\neq j$ with $U_i=U_j$:} Then $\a_i\neq\a_j$; suppose $\a_i>\a_j$. Hence (\ref{43}) implies that (writing $U$ for $U_i=U_j$)
$x_i\in S_{\a_i,m}(U)\cap V_{\a_j}(U)$, contradicting condition (b${}_U$).

\smallskip

\noindent{\em Case 2.}  {\em The $U_j$ are all distinct:} Let $E=\{x_j\st j\in\NNNN\}$. Since $x_i\in V_{\a_j}(U_j)\sbe U_j$ for all $i,j\in\NNNN$, we have that $E\sbe\bigcap\{U_j\st j\in\NNNN\}$. But each $U_j$ belongs to $\UU_n$ and $\UU_n$ is point-finite. Hence, since the $U_j$ are all distinct, $\bigcap\{U_j\st j\in\NNNN\}=\ems$. Since $E\nes$, we get a contradiction.

\smallskip

Therefore, all $\WW_{m,n,k}$ are point-finite.
\sqs

A modification of the proof of the Michael-Rudin Theorem
\ref{MR1Th14} yields the next result:

\begin{theorem}\label{th87}{\rm (\cite{D0})} A compact Hausdorff space $X$ is a  UEC iff it has an F-separating
$\sigma$-boundedly point-finite family of open subsets.
\end{theorem}

\doc ($\Rightarrow$) By  Theorem \ref{argf}, $X$ has a $T_0$-separating $\s$-boundedly point-finite family of cozero subsets. Now, exactly as in the proof of
the necessity of
Theorem \ref{MR1Th14}, we get that $X$ has an F-separating $\s$-boundedly point-finite family of open subsets.

\smallskip

($\Leftarrow$)
Let $\UU=\bigcup\{\UU_n\st n\in\NNNN\}$ be an F-separating open family in $X$, where all $\UU_n$ are boundedly point-finite collections.
We may suppose that $\UU_n\sbe \UU_{n+1}$ and $\UU_n$ covers $X$, for all $n\in\NNNN $.

For every $n\in\NNNN$, let $$\bn=1+ord(\UU_n).$$

Now, we define, as in the proof of Theorem \ref{MR1Th14}, the family $\WW$ and the families $\WW_{m,n,k}$, where $m,n,k\in\NNNN$. Then, as it it shown
in the proof of Theorem \ref{MR1Th14}, $\WW$ is a $T_0$-separating family in $X$, $\WW=\bigcup\{\WW_{m,n,k}\st m,n,k\in\NNNN\}$ and $\WW$ consists of cozero subsets of $X$. We will prove that all families $\WW_{m,n,k}$ are boundedly point-finite. In what follows, we shall use the notation from the proof of Theorem \ref{MR1Th14}.

We shall prove that for all $m,n,k\in\NNNN$,
$$ord(\WW_{m,n,k}) < 2^{\bk}\cdot\bn.$$

Let $m,n,k\in\NNNN$. Suppose that there exists $x^*\in X$ such that $|\WW_{m,n,k})(x^*)|> 2^{\bk}\cdot\bn$. Let $n_0=1+2^{\bk}\cdot\bn$. Then there are distinct
$W_{\a_j,m}(U_j)$, $j=1,\ldots, n_0$, in $\WW_{m,n,k}$ all containing the point $x^*$. Set $J=\{1,\ldots,n_0\}$. Then $x^*\in H_{\a_j,m}(U_j)$ for every $j\in J$. Hence there are points $x_j\in S_{\a_j,m}(U_j)$, $j\in J$, such that $x^*\in\bigcap\UU_k(x_j)$ for all $j\in J$. Therefore, $\UU_k(x^*)\spe \UU_k(x_j)$ for all $j\in J$.
We have that $|\UU_k(x^*)|\le\bk$. Let $P=exp(\UU_k(x^*))$. Then, setting, for all $j\in J$, $\p(j)=\UU_k(x_j)$, we define a function
$$\p:J\lra P.$$
Suppose that $|\p\inv(p)|\le\bn$ for all $p\in P$. Then we will obtain that $n_0=|J|\le|P|\cdot\bn\le 2^{\bk}\cdot \bn< n_0$, which is a contradiction.
  Hence, there exists
$p_0\in P$ such that $|\p\inv(p_0)|>\bn$. This implies that there exists a subset $I$ of $J$ such that $|I|=\bn + 1$ and
$\UU_k(x_i)=\UU_k(x_j)$ for all $i,j\in I$. Then, by (\ref{42}), we obtain that $x_i\in V_{\a_j}(U_j)$ for all $i,j\in I$. Hence, we have that
\begin{equation}\label{433}
x_i\in S_{\a_i,m}(U_i)\cap V_{\a_j}(U_j) \mbox{ for all } i,j\in I.
\end{equation}

Now there are two cases.

\smallskip

\noindent{\em Case 1.}  {\em There are $i\neq j$, $i,j\in I$, with $U_i=U_j$:} Then $\a_i\neq\a_j$
(because all the sets $W_{a_l,m}(U_l)$, $l\in J$, are distinct). Suppose $\a_i>\a_j$. Then (\ref{433}) implies that
(writing $U$ for $U_i=U_j$) $x_i\in S_{\a_i,m}(U)\cap V_{\a_j}(U)$, contradicting condition (b${}_U$)  from the proof of Theorem \ref{MR1Th14}.

\smallskip

\noindent{\em Case 2.}  {\em The $U_j$ are all distinct, $j\in I$:} Let $E=\{x_j\st j\in I\}$. Since $x_i\in V_{\a_j}(U_j)$ for all $i,j\in I$ (by (\ref{433})), we have that $E\sbe\bigcap\{U_j\st j\in I\}$. However all $U_j$ belong to $\UU_n$ and $ord(\UU_n)=\bn - 1$. Hence $\bigcap\{U_j\st j\in I\}=\ems$. Since $E\nes$, we get a contradiction.

Therefore, all families $\WW_{m,n,k}$ are boundedly point-finite. Thus $\WW$ is a $T_0$-separating $\s$-boundedly point-finite family of cozero subsets of $X$.
Now, Theorem \ref{argf} implies that $X$ is a UEC.
\sqs

Following the ideas of the proof of Theorem \ref{MR2Th11} of E. Michael and M. E. Rudin, we obtain
such a generalization of it:

\begin{theorem}\label{GMR}{\rm (\cite[Theorem 3]{D1})}
If $X$ is compact Hausdorff, and if $X=X_1 \cup X_2$ with $X_1$ and $X_2$ having
$\s$-strongly point-finite (resp. $\s$-boundedly point-finite) bases, then $\cl(X_1)\cap\cl(X_2)$ is metrizable and
$X$ is an EC (resp., a UEC).
\end{theorem}

\doc For $i=1,2$, let $\UU_i$ be a $\s$-strongly point-finite (resp., $\s$-boundedly point-finite) base for $X_i$. For $i=1,2$ and any $U\in\UU_i$, let $\Phi_i(U)$ be an open subset of $X$ such that $\Phi_i(U)\cap X_i=U$. Let
$$\UU=\{\Phi_i(U)\cap\cl(X_1)\cap\cl(X_2)\st U\in\UU_i, i=1,2\}.$$
We will prove that $\UU$ is a $\s$-point-finite (resp. $\s$-boundedly point-finite) base of $X_0=\cl(X_1)\cap\cl(X_2)$.

Let us first show that $\UU$ is a base of $X_0$. Since $\cl(X_1)\stm X_1\sbe X_2$, we have that $(\cl(X_1)\stm X_1)\cap(\cl(X_2)\stm X_2)=\ems$. Thus
$X_0=(X_1\cap\cl(X_2))\cup(\cl(X_1)\cap X_2)$. Let $x\in X_0$, $O_o$ be  an open neighborhood of $x$ in $X_0$ and $O$ be an open subset of $X$ such that $O\cap X_0=O_0$. Let $V$ be  an open neighborhood of $x$ in $X$ such that $\cl(V)\sbe O$. Now there are two possibilities:

\smallskip

\noindent{\it Case 1.} $x\in X_1\cap\cl(X_2)$. Then there exists $U\in\UU_1$ such that $x\in U\sbe V\cap X_1$. We will show that $x\in \Phi_1(U)\cap X_0\sbe O_0$.
Obviously, $x\in\Phi_1(U)\cap X_0$. Further, $\Phi_1(U)\cap \cl(X_1)\sbe Ex_{cl(X_1)}(U)\sbe \cl_{cl(X_1)}(U)=\cl_X(U)$ and thus
$\Phi_1(U)\cap X_0\sbe \cl(U)\cap\cl(X_2)\sbe  \cl(V)\cap\cl(X_2)\sbe \cl(V)\sbe O$. Hence, $x\in \Phi_1(U)\cap X_0\sbe O\cap X_0=O_0$ and
$\Phi_1(U)\cap X_0\in\UU$.

\smallskip

\noindent{\it Case 2.} $x\in X_2\cap\cl(X_1)$. Then the proof is analogous to that of Case 1.

\smallskip

So, $\UU$ is a base of $X_0$. Now,
Lemma \ref{Ex}
implies that for $i=1,2$, $Ex_{cl(X_i)}(\UU_i)$ is a $\s$-point-finite family (resp., $\s$-boundedly point-finite family) and, hence, $\UU$ is such a family too. Therefore, $X_0=\cl(X_1)\cap\cl(X_2)$ has a point-countable base. Then, by the Mi\v{s}\v{c}enko Theorem (see, e.g., \cite[3.12.23(f)]{E}), $X_0$ is a metrizable space.

Let $\BB$ be a countable base of $X_0$. Then for every $V_1, V_2\in\BB$ such that $\cl_{X_0}(V_1)\cap \cl_{X_0}(V_2)=\ems$ (equivalently, $\cl_{X}(V_1)\cap \cl_{X}(V_2)=\ems$) there exists an open subset $W_{V_1,V_2}$ of $X$ such that $\cl_X(V_1)\sbe W_{V_1,V_2}\sbe \cl_X(W_{V_1,V_2})\sbe X\stm  \cl_X(V_2)$. Let $\VV=\{W_{V_1,V_2}\st V_1, V_2\in\BB, \cl_{X}(V_1)\cap \cl_{X}(V_2)=\ems\}$. Let $Y=X\stm X_0$. Then $Y$ is an open subset of $X$ and $Y$ is a disjoint sum of the spaces $X_1\cap Y$ and $X_2\cap Y$. Since the spaces $X_1\cap Y$ and $X_2\cap Y$ have $\s$-point-finite (resp. $\s$-boundedly point-finite) bases, the space $Y$ also has such a base $\WW$. Moreover, we can suppose that $\cl_X(W)\sbe Y$, for every $W\in\WW$. Let $\Omega=\VV\cup\WW$. Then, obviously, $\Omega$ is an F-separating open family in $X$. Moreover, it is a
$\s$-point-finite (resp. $\s$-boundedly point-finite) family. Therefore, by the Michael-Rudin Theorem
\ref{MR1Th14}
(resp., by Theorem \ref{th87}), we get that $X$ is an EC (resp., UEC).
\sqs

\begin{cor}\label{th84} If $X$ has a $\sigma$-strongly point-finite (resp., $\sigma$-boundedly point-finite)  base and has a
compactification $bX$ with remainder $Y$ which also has a
$\sigma$-strongly point-finite (resp., $\sigma$-boundedly point-finite) base, then $bX$ is an EC (resp., UEC) and $\cl_{bX}(Y)$
is metrizable.
\end{cor}

\section{On a theorem of D. Preiss and P. Simon}

 We are now going to generalize the following theorem of D. Preiss and P. Simon:

\begin{theorem}\label{e6}{\rm (\cite[Corollary 7]{PrS})}
 If $X$ is EC and $Y$ is a dense proper subspace of $X$, then $X\neq \b Y$.
 \end{theorem}

It follows immediately from another theorem of D. Preiss and P. Simon:

\begin{theorem}\label{e66}{\rm (\cite[Corollary 6]{PrS})}
 A pseudocompact subspace of an Eberlein compact is closed and hence it is an Eberlein compact, too.
 \end{theorem}

Indeed, as it is stated in \cite{PrS}, for proving Theorem \ref{e6} one can argue as follows: if $X=\b Y$ then, obviously, $X=\b(X\stm \{x\})$ for any $x\in X\stm Y$; then $X\stm\{x\}$ is a non-closed pseudocompact subspace (see, e.g., \cite[3.12.16(b)]{E}) of $X$, which contradicts Theorem \ref{e66}.

\begin{rem}\label{remps}
\rm
Let us note that Theorem \ref{e66} can be also proved  using the Yakovlev Theorem (\cite[Corollary 2]{Y}) (= every Eberlein compact is hereditarily $\s$-metacom\-pact)
and Uspenskii's Theorem (\cite[Proposition]{Usp}) (= every $\s$-point-finite open cover of a pseudocompact space $X$ contains a finite subfamily whose union is dense in $X$)  (the later implies immediately that every $\s$-metacompact pseudocompact regular space is compact). Note that the cited theorems of Yakovlev and Uspenskii were published many years after the publication of the Preiss-Simon Theorem \ref{e66}.
\end{rem}

We will need the following definition from \cite{D1}:

\begin{defi}\label{femb}{\rm (\cite[Definition 5]{D1})}
\rm
 Let $X$ be a space, $Y\sbe X$ and $\tau$ be a cardinal number. The subspace $Y$ is said to be
{\em $F_\tau$-embedded in $X$}\/ if there exists a family $\UU$ which is the union of $\tau$ point-finite open
families in $X$ and which F-separates $Y$.
\end{defi}

Now, we will prove a theorem which, as we will see below in Remark \ref{rempscor}, generalizes the Preiss-Simon Theorem \ref{e6}.

\begin{theorem}\label{GPS0}{\rm (\cite[Theorem 4]{D1})}
 Let $X$ be compact Hausdorff and  $X=Y \cup Z$, where $Y$ is a  proper dense E-subspace
of $X$ and  $Z$ is $F_{\aleph_0}$-embedded in $X$. Then $X\neq \b Y$.
\end{theorem}

\doc  By Definition \ref{espa} and Remark \ref{rem0}, there exists a $\s$-strongly point-finite F-separating almost subbase $\a_Y$ of $Y$. Since $Z$ is $F_{\aleph_0}$-embedded in $X$, there exists a $\s$-point-finite family $\a_Z$ of open subsets of $X$ which F-separates $Z$.

Suppose that $X=\b Y$. Let
$$\a=\{Ex_X(V)\st V\in\a_Y\}\cup\a_Z.$$
Note that, by Lemma \ref{Ex}, 
\begin{equation}\label{4ps}
\mbox{the family } Ex_X(\a_Y) \mbox{ is } \s\mbox{-point-finite} 
\end{equation}
and thus the family $\a$ is $\s$-point-finite.
For every $V\in\a_Y$ there exists a U-representation $U(V)=\{U_n(V)\st n\in\NNNN\}$ of $V$ in $Y$ such that the family
$$\a_Y\ap=\a_Y\cup\{Y\stm U_{2n-1}(V)\st V\in\a_Y,n\in\NNNN\}$$
is a subbase of $Y$; for every $n\in\NNNN$, we have that $U_{2n-1}(V)$ is a zero-set in $Y$, $U_{2n}(V)$ is a cozero-set in $Y$ and $U_{2n-1}(V)\sbe U_{2n}(V)$. Then, by \cite[Theorem 1.5.14]{E}, for every $V\in\a_Y$ and for every $n\in\NNNN$ there exists a continuous function $f_{n,V}:Y\lra \IIII$ such that $f_{n,V}\inv(0)=U_{2n-1}(V)$ and $f_{n,V}\inv(1)=Y\stm U_{2n}(V)$.
Set, for every $V\in\a_Y$ and every $n,k\in\NNNN$, $U_{2n-1}^{2k-1}(V)=f_{n,V}\inv([0,\frac{1}{2k}])$ and $U_{2n-1}^{2k}(V)=f_{n,V}\inv([0,\frac{1}{2k}))$.
 Then
 $U_{2n-1}(V)=\bigcap\{U_{2n-1}^k(V)\st k\in\NNNN\}$, 
 \begin{equation}\label{5ps}
 U_{2n-1}^{k+1}(V)\sbe U_{2n-1}^{k}(V),\ \  U_{2n-1}^{1}(V)\sbe U_{2n}(V), 
 \end{equation}
 $U_{2n-1}^{2k-1}(V)$ is a zero-set in $Y$ and $U_{2n-1}^{2k}(V)$ is a cozero-set in $Y$, for every $k\in\NNNN$.
Let
$$\a\ap=\a\cup\{Ex_X U_{2n}(V)\st V\in\a_Y,n\in\NNNN\}\cup\{Ex_X U_{2n-1}^{2k}(V)\st V\in\a_Y,n,k\in\NNNN\}.$$
We will show that $\a\ap$ is an F-separating family in $X$. Indeed, let $x,y\in X$ and $x\neq y$.

\smallskip

\noindent{\it Case 1:} $x,y\in Y$.
Then there exists $V\in\a_Y$ such that $x\in V$ and $y\nin\cl_Y(V)$ or vice versa. Thus $x\in Ex_X(V)\in\a\sbe\a\ap$ and $y\nin\cl_X(V)=\cl_X(Ex_X(V))$ or vice versa.

\smallskip

\noindent{\it Case 2:} $x,y\in Z$.
Then there exists $W\in\a_Z\sbe\a\ap$ such that $x\in W$, $y\nin\cl_X(W)$ or vice versa.

\smallskip

\noindent{\it Case 3:} $x\in Y$, $y\in Z$ {\em or vice versa.} Clearly, it is enough to consider the case when
 $x\in Y$ and $y\in Z$. There exist disjoint open in $X$ neighborhoods $Ox$ and $Oy$ respectively of $x$ and $y$. Now, there exists a finite subfamily $\b_x$ of the family $\a_Y\ap$ such that $x\in\bigcap\b_x\sbe Y\cap Ox$. Then $y\nin\cl_X(\bigcap\b_x)$. Let $\b_x=\{V_i\in\a_Y\st i=1,\ldots, k\}\cup \{Y\stm U_{2n_j -1}(W_j)\st W_j\in\a_Y,j=1,\ldots, s\}.$ Then for every $i\in\{1,\ldots,k\}$, there exists $m_i\in\NNNN$ such that $x\in U_{2m_i}(V_i)$; also, for every $j\in\{1,\ldots,s\}$, there exists $l_j\in\NNNN$ such that $x\in Y\stm U_{2n_j-1}^{2l_j-1}(W_j)$. Then $x\in\bigcap\{U_{2m_i}(V_i)\st i=1,\ldots, k\}\cap\bigcap\{Y\stm U_{2n_j-1}^{2l_j-1}(W_j)\st j=1,\ldots, s\}\sbe
\bigcap\{U_{2m_i+1}(V_i)\st i=1,\ldots, k\}\cap\bigcap\{Y\stm U_{2n_j-1}^{2l_j}(W_j)\st j=1,\ldots, s\}\sbe\bigcap\b_x\sbe Ox.$
Hence, since $X=\b Y$ (and, thus, we can use \cite[Theorem 1.46(4)]{Walker}), we obtain that $y\nin\cl_X(\bigcap\{U_{2m_i+1}(V_i)\st i=1,\ldots, k\}\cap\bigcap\{Y\stm U_{2n_j-1}^{2l_j}(W_j)\st j=1,\ldots, s\})=
\bigcap\{\cl_X(U_{2m_i+1}(V_i))\st i=1,\ldots, k\}\cap\bigcap\{\cl_X(Y\stm U_{2n_j-1}^{2l_j}(W_j))\st j=1,\ldots, s\}$. Thus we have the following two possibilities:

\smallskip

\noindent(a) {\em There exists $i_0\in\{1,\ldots,k\}$ such that} $y\nin \cl_X(U_{2m_{i_0}+1}(V_{i_0})).$

Then we have that $x\in Ex_X(U_{2m_{i_0}}(V_{i_0}))\in\a\ap$, $$\cl_X(Ex_X(U_{2m_{i_0}}(V_{i_0})))=\cl_X(U_{2m_{i_0}}(V_{i_0}))\sbe\cl_X(U_{2m_{i_0}+1}(V_{i_0}))$$ and hence $y\nin\cl_X(Ex_X(U_{2m_{i_0}}(V_{i_0})))$.

\smallskip

\noindent(b) {\em There exists $j_0\in\{1,\ldots,s\}$ such that} $y\nin \cl_X(Y\stm U_{2n_{j_0}-1}^{2l_{j_0}}(W_{j_0})).$

Then $y\in Ex_X(U_{2n_{j_0}-1}^{2l_{j_0}}(W_{j_0}))\in\a\ap$. Since $x\nin U_{2n_{j_0}-1}^{2l_{j_0}-1}(W_{j_0})$ and
$$U_{2n_{j_0}-1}^{2l_{j_0}}(W_{j_0})\sbe U_{2n_{j_0}-1}^{2l_{j_0}-1}(W_{j_0}),$$
we obtain that $x\nin\cl_Y(U_{2n_{j_0}-1}^{2l_{j_0}}(W_{j_0}))$ and thus $x\nin\cl_X(Ex_X(U_{2n_{j_0}-1}^{2l_{j_0}}(W_{j_0})))$.

Therefore $\a\ap$ is an F-separating family of open subsets of $X$. It is easy to see, using  (\ref{4ps}), (\ref{5ps}) and Definition \ref{Urepr}, that $\a\ap$ is also a $\s$-point-finite family in $X$. Thus, by Theorem
\ref{MR1Th14}, $X$ is an Eberlein compact. Since $X=\b Y$ and $Y\neq X$, using Theorem \ref{e6}, we get a contradiction.
\sqs

\begin{rem}\label{rempscor}
\rm
Let us note that Theorem \ref{GPS0} is a generalization of Theorem \ref{e6}. Indeed, let $X$ be an EC and $Y$ be a dense proper subspace of $X$. Suppose that $X=\b Y$. Since $Y\neq X$, there exists a point $x\in X\stm Y$. Set $Z=\{x\}$ and $Y\ap=X\stm\{x\}$. Then, by Corollary \ref{subsp}, $Y\ap$ is an E-space. Obviously, $Z$ is $F_{\aleph_0}$-embedded in $X$, $X=Y\ap\cup Z$ and $X=\b Y\ap$. Thus, by Theorem \ref{GPS0}, we get a contradiction. So, $X\neq\b Y$.
\end{rem}

\begin{rem}\label{rempscor1}
\rm
Clearly, if $Y$ is $F_\tau$-embedded in $X$ and $Z\sbe Y$ then $Z$ is $F_\tau$-embedded in $X$. Thus, Theorem \ref{GPS0} can be reformulated as follows:

{\em Let $X$ be a non-compact E-space, $cX$ be a Hausdorff compactification of $X$ and $cX\stm X$ be $F_{\aleph_0}$-embedded in $cX$. Then $cX\neq\b X$.}
\end{rem}

Now we will list some sufficient conditions which a subspace $Y$ of a space $X$ has to satisfy in order to be $F_\tau$-embedded in $X$, where $\tau$ is a cardinal number.

Recall that if $X$ is a topological space then by $nw(X)$ is denoted the network weight of $X$ (see, e.g., \cite{E}).

\begin{theorem}\label{Ftauemb}{\rm (\cite{D0})} Let $X$ be a space, $Y$ be a subspace of $X$ and $\tau$ be a cardinal number. Then $Y$ is $F_\tau$-embedded in $X$ if some of the following conditions holds:

\smallskip

\noindent(a) $X$ is a normal $T_1$-space, $\tau\ge\aleph_0$ and $nw(Y)\le\tau$;

\smallskip

\noindent(b)  $Y$ is a preopen subset of $X$ (i.e.
$Y\sbe \int_X(\cl_X(Y))$)  and there is a family $\a$ of open in $Y$ subsets of $Y$ which is the union of $\tau$ strongly point-finite collections and which F-separates   $Y$.
\end{theorem}

\doc (a) Let $\NN_Y$ be a network in $Y$ with $|\NN_Y|\le\tau$ and let $\OM=\{(A,B)\st A,B\in\NN_Y,\cl_X(A)\cap\cl_X(B)=\ems\}$. For every pair $(A,B)\in\OM$ we fix an open subset $W_{(A,B)}$ of $X$ such that $\cl_X(A)\sbe W_{(A,B)}\sbe\cl_X(W_{(A,B)})\sbe X\stm \cl_X(B)$. Then $\WW=\{W_{(A,B)}\st (A,B)\in\OM\}$ is a collection of open subsets of $X$ and $|\WW|\le\tau$ (thus, $\WW$ is the union of $\tau$  point-finite collections of open subsets of $X$). We will show that $\WW$ F-separates $Y$. Indeed, let $x,y\in Y$ and $x\neq y$. Then there are two  open neighborhoods $Ox$ and $Oy$ respectively of $x$ and $y$ such that $\cl_X(Ox)\cap\cl_X(Oy)=\ems$. There exist $A,B\in\NN_Y$ such that $x\in A\sbe Y\cap Ox$ and $y\in B\sbe Y\cap Oy$. Then $(A,B)\in\OM$, $x\in W_{(A,B)}$ and $y\nin\cl_X(W_{(A,B)})$. Therefore, $Y$ is $F_\tau$-embedded in $X$.

\smallskip

\noindent(b)  Let $Z=\int_X(\cl_X(Y))$. Then $Y\sbe Z$, $Y$ is dense in $Z$ and $Z$ is open in $X$. Then, by Lemma \ref{Ex}, the family $\a\ap=\{Ex_Z(U)\st U\in\a\}$ is the union of $\tau$ point-finite collections of open subsets of $X$. Obviously, the family $\a\ap$ F-separates $Y$. Therefore,  $Y$ is $F_\tau$-embedded in $X$.
\sqs

\begin{cor}\label{e55}{\rm (\cite{D0})}
Let $X$ be a space and $Y$ be a preopen (e.g., dense) E-subspace of $X$. Then $Y$ is $F_{\aleph_0}$-embedded in $X$. In particular, every metrizable dense subspace $Y$ of $X$ is $F_{\aleph_0}$-embedded in $X$.
\end{cor}

\doc It follows from Definition \ref{espa}, Remark \ref{rem0}  and Theorem \ref{Ftauemb}(b).
\sqs

\begin{cor}\label{e5}{\rm (\cite[Corollary 8]{D1} and \cite{D0})}
 Let $X$ be a compact Hausdorff space and $X=Y \cup Z$, where $Y$ is a proper dense E-subspace
of $X$. Then $X\neq\b Y$  if some of the following conditions holds:

\smallskip

\noindent(a) $nw(Z)\le\aleph_0$;

\smallskip

\noindent(b) $w(Z)\le\aleph_0$;

\smallskip

\noindent(c) $Z$ is preopen in $X$ and there is a
$\s$-strongly point-finite  family of open subsets of $Z$ which F-separates $Z$;

\smallskip

\noindent(d) $Z$ is a preopen (e.g., dense) E-subspace of $X$;

\end{cor}

\doc All follows from Theorem \ref{GPS0}, Theorem \ref{Ftauemb} and Corollary \ref{e55}.
\sqs

\begin{cor}\label{e567}
Let $X$ be a non-compact E-space and $X=R(X)$. Then $\b X\stm X$ is not an E-space.
\end{cor}

\doc  Suppose that $\b X\stm X$ is an E-space. Since $X=R(X)$, we get that $\cl_{\b X}(\b X\stm X)=\b X$ (see, e.g., \cite{AP}),
i.e. $\b X\stm X$ is a dense subspace of $\b X$.
Thus, by Corollary \ref{e5}(d), $X$ is not C*-embedded in $\b X$ - a contradiction. Therefore $\b X\stm X$ is not an E-space.
\sqs

Let us recall a result of P. Simon which will be used below:

\begin{theorem}\label{ps-sec}{\rm (\cite[Proposition 9]{Sim})}
A compact Hausdorff space is a SEC if and only if it has a point-finite $T_0$-separating
collection of open sets.
\end{theorem}

Recall that a subset $S$ of a topological space $X$ is {\em zero-dimensionally embedded in} $X$ if $X$ has a base $\BB$ for the open sets such that $Fr_X(U)\sbe X\stm S$ for all $U\in\BB$.
Also, we have to recall the well-known fact that if $X$ is a {\em rim-compact space} (= {\em semi-compact space} = {\em peripherically-compact space}\/) (i.e., a Hausdorff space which has a base of open sets with compact boundaries) and $FX$ is its {\em Freudenthal compactification}\/ (\cite{Fre,Mo,AaNi}) (i.e., the largest compactification of $X$ with a zero-dimensionally embedded remainder)  then any two disjoint closed subsets of $X$ with compact boundaries have disjoint closures in $FX$; moreover, if $F$ and $G$ are closed subsets of $X$ with compact boundaries then $\cl_{FX}(F\cap G)=\cl_{FX}(F)\cap\cl_{FX}(G)$.

The proof of the next theorem is similar to that of Theorem \ref{GPS0}.

\begin{theorem}\label{fr}{\rm (\cite[Theorem 5]{D1})}
Let $X$ be a non-compact SE-space, $FX$ be its Freudenthal compactification and $FX =X \cup Y$,
where $Y$ is $F_1$-embedded (resp., $F_{\aleph_0}$-embedded) in $FX$. Then $FX$ is a SEC (resp., EC).
\end{theorem}

\doc There exists a strongly point-finite family $\a_X$ of clopen subsets of $X$ such that the family
$$\a\ap_X=\a_X\cup\{X\stm V\st V\in\a_X\}$$
is a subbase of $X$.

\smallskip

\noindent(A) {\em Let $Y$ be $F_1$-embedded in $FX$.} Then there exists a point-finite family $\a_Y$ of open subsets of $FX$ which F-separates $Y$. Set
$$\a=\a_Y\cup\{Ex_{FX}(V)\st V\in\a_X\}.$$
We will show that $\a$ is an F-separating family in $FX$. Indeed, let $x,y\in FX$ and $x\neq y$. There are three cases.

\smallskip

\noindent{\em Case 1: $x,y\in X$.} Then there exists a finite subfamily $\b$ of the family $\a\ap_X$ such that $x\in\bigcap\b\sbe X\stm\{y\}$. Obviously, there exists $V\in\b$ such that $y\nin V$. If $V\in\a_X$ then we have that $x\in Ex_{FX}(V)$, $Ex_{FX}(V)\in\a$ and $y\nin\cl_{FX}(Ex_{FX}(V))$.
If $V\nin\a_X$ then $X\stm V\in\a_X$. Thus, setting $W=X\stm V$, we get that
$y\in Ex_{FX}(W)$, $Ex_{FX}(W)\in\a$ and $x\nin\cl_{FX}(Ex_{FX}(W))$.

\smallskip

\noindent{\em Case 2: $x,y\in Y$.} Then there exists $V\in\a_Y\sbe\a$ such that $x\in V$ and $y\nin\cl_{FX}(V)$ or vice versa.

\smallskip

\noindent{\em Case 3: $x\in X$ and $y\in Y$ or vice versa.} Clearly, it is enough to regard  the case when  $x\in X$ and $y\in Y$. There exists disjoint open in $FX$ neighborhoods $Ox$ and $Oy$ of $x$ and $y$, respectively. There exists a finite subfamily $\b$ of the family $\a\ap_X$ such that $x\in\bigcap\b\sbe X\cap Ox$. Then $y\nin\cl_{FX}(\bigcap\b)$. Let $\b=\{V_i\in\a_X\st i=1,\ldots,k\}\cup\{X\stm W_j\st W_j\in\a_X, j=1,\ldots,l\}$. Then $\cl_{FX}(\bigcap\b)=\bigcap\{\cl_{FX}(V_i)\st i=1,\ldots,k\}\cap\bigcap\{\cl_{FX}(X\stm W_j)\st j=1,\ldots,l\}$. There are two possibilities.

\smallskip

\noindent(a) {\em There exists $i\in\{1,\ldots,k\}$ such that $y\nin\cl_{FX}(V_i)$.} Then $x\in V_i\sbe Ex_{FX}(V_i)\in\a$ and $y\nin\cl_{FX}(V_i)=\cl_{FX}(Ex_{FX}(V_i))$.

\smallskip

\noindent(b) {\em There exists $j\in\{1,\ldots,l\}$ such that $y\nin\cl_{FX}(X\stm W_j)$.} Then $y\in FX\stm cl_{FX}(X\stm W_j)=Ex_{FX}(W_j)\in\a$. Since $x\nin W_j$ and $W_j$ is closed in $X$, we get that $x\nin\cl_{FX}(W_j)$. Note, finally, that $\cl_{FX}(W_j)=\cl_{FX}(Ex_{FX}(W_j))$. Therefore, $x$ and $y$ are F-separated by $\a$.

\smallskip

\smallskip

So, the family $\a$ F-separates $FX$. By Lemma \ref{Ex}, the family $Ex_{FX}(\a_X)$ is point-finite. Hence, the family $\a=\a_Y\cup Ex_{FX}(\a_X)$ is point-finite. Now, Theorem \ref{ps-sec} implies that $FX$ is a SEC.

\smallskip

\noindent(B) {\em Let $Y$ be $F_{\aleph_0}$-embedded in $FX$.} Then there exists a $\s$-point-finite family $\a_Y$ of open subsets of $FX$ which F-separates $Y$. Now, we show exactly as in (A) that the family
$$\a=\a_Y\cup Ex_{FX}(\a_X)$$
F-separates $FX$. Using again Lemma \ref{Ex}, we get that $\a$ is a $\s$-point-finite family of open subsets of $FX$. Thus, by the Michael-Rudin Theorem \ref{MR1Th14}, we get that $FX$ is an EC.
\sqs

\begin{rem}\label{rfr1}
\rm
We can reformulate Theorem \ref{fr} as follows:

{\em Let $X$ be a non-compact SE-space, $FX$ be its Freudenthal compactification and $FX \stm X$ be
 $F_1$-embedded (resp., $F_{\aleph_0}$-embedded) in $FX$. Then $FX$ is a SEC (resp., EC).}
\end{rem}

\begin{cor}\label{e7}{\rm (\cite{D0})}
If $X$ is an an SE-space and $nw(FX\stm X)\le\aleph_0$, then $FX$ is EC.
\end{cor}

\doc It follows from Theorem \ref{fr} and Theorem \ref{Ftauemb}(a).
\sqs

\begin{cor}\label{e77}{\rm (\cite[Corollary 10]{D1})}
If $X$ is  SE-space and $w(FX\stm X)\le\aleph_0$, then $FX$ is EC.
\end{cor}

\section{On a question of A. V. Arhangel'ski\v{\i}}

During the his talk on the international conference $``$Analysis, Topology and Applications - 2014", Vrnja\v{c}ka Banja, Serbia, May 25-29, 2014, A. V. Arhangel'ski\v{\i} asked the following question:

\begin{question}\label{p2014}
\rm
Is it true that if a Tychonoff space $X$ has a $\s$-disjoint base then it can be embedded in an Eberlein compact?
\end{question}

Let us recall that every perfectly normal space with a $\sigma$-disjoint base is metrizable (see \cite{Aull}) and thus, by the Arhangel'ski\v{\i} Theorem \cite[Theorem 14]{A}, it has a compactification which is an EC.

The next result follows immediately from Theorem \ref{MainE}:

\begin{theorem}\label{c87} Let $X$ be a Tychonoff space. Then:

\smallskip

 \noindent(a) If $X$ has a $\sigma$-disjoint  base of cozero-sets then $X$ can be embedded in a uniform Eberlein compact;

 \smallskip

 \noindent(b) If  $X$ has a $\sigma$-strongly point-finite  base of cozero-sets then $X$ can be embedded in an Eberlein compact.
\end{theorem}

With the next example we show that the class of Tychonoff spaces with a $\sigma$-disjoint base of cozero-sets is strictly larger than the class of metrizable spaces.

 \begin{exa}\label{ex81}
 \rm
 There exists a non-metrizable Tychonoff space $X$ which has a  $\sigma$-disjoint base consisting of cozero-sets.
Such are, e.g., the following spaces:

\smallskip

\noindent(a) The Michael line $\RRRR_\QQQQ$;

\smallskip

\noindent(b) The Corson-Michael example $(X,\tau)$ \cite[Example 6.4]{CM} of a non-metrizable, hereditarily paracompact, Lindel\"{o}f space with a $\s$-disjoint base.
\end{exa}

\doc  (a) For a description of the Michael line $\mathbb{R}_\mathbb{Q}$ see, e.g., \cite[Example 5.1.22]{E}.
Let $\BB\ap$ be a $\s$-discrete base of $\mathbb{R}$ with its natural topology and $\PP=\{\{x\}\st x\in\mathbb{P}\}$. Then $\BB=\BB\ap\cup\PP$ is a $\sigma$-disjoint base of the Michael line $\mathbb{R}_\mathbb{Q}$ consisting of open $F_\s$-subsets of $\RRRR_\QQQQ$ and thus, since $\RRRR_\QQQQ$ is a normal space, of cozero-sets. Finally, by \cite[Problem 5.5.2(d)]{E}, $\RRRR_\QQQQ$ is a non-metrizable space.

\smallskip

\noindent(b)  We will use the notation from \cite[Example 6.4]{CM}. The required  $\sigma$-disjoint base consisting of cozero-sets can be constructed exactly as in (a):
let $\BB\ap$ be a $\s$-discrete base of $(X,\s)$  and $\PP=\{\{y\}\st y\in Y\}$. Then $\BB=\BB\ap\cup\PP$ is a $\sigma$-disjoint base of $(X,\tau)$ consisting of open $F_\s$-subsets of $(X,\tau)$ and thus, since $(X,\tau)$ is a normal space, of cozero-sets.
\sqs

Another partial answer to the Problem \ref{p2014} of A. V. Arhangel'ski\v{\i} gives the following theorem.

\begin{theorem}\label{th83} Every  space $X$ which has a $\sigma$-strongly point-finite family
(resp., $\sigma$-boundedly point-finite family)
$\PP$ consisting of clopen sets, such that $\PP\cup \{X\setminus U: U\in \PP\}$ is a subbase of $X$, can be embedded in a zero-dimensional  Eberlein compact (resp., zero-dimensional uniform Eberlein compact). In particular, every  space $X$ which has a $\sigma$-strongly point-finite
(resp., $\sigma$-boundedly point-finite;  $\s$-disjoint)
base consisting of clopen sets can be embedded in a zero-dimensional  Eberlein compact (resp., zero-dimensional uniform Eberlein compact).
\end{theorem}

\doc Let $\BB$ be the Boolean subalgebra of $CO(X)$ generated by the $\sigma$-strongly point-finite family $\PP$. Then, clearly, $\BB$ is a base of $X$. Let $Y=S(\BB)$, i.e., $Y$ is the Stone space of the Boolean algebra $\BB$. Then $Y$ is a compactification of $X$ (see \cite[Theorem 13.1 and the paragraph before it]{Dw}). We shall show that $Y$ is an Eberlein compact. By the Stone Duality Theorem \cite{ST}, there exists a Boolean isomorphism $\p:\BB\lra CO(Y)$. Let $\QQ=\p(\PP)$, $\QQ^*=\{Y\stm Q\st Q\in\QQ\}$ and $\PP^*=\{X\setminus U: U\in \PP\}$. Obviously, $\QQ$ is a $\s$-strongly point-finite family. Thus, having in mind Theorem \ref{unifthE}, we need only to show that $\QQ$ is an almost subbase of $Y$. Clearly, if $V\in\QQ$ and we set $U_n(V)=V$ for every $n\in\NNNN$, then $\{U_n(V)\st n\in\NNNN\}$ is a U-representation of $V$. Hence, it is enough to show that the family $\QQ\cup
\QQ^*$ is a subbase of $Y$. Obviously, $\p(\PP\cup\PP^*)=\QQ\cup
\QQ^*$. Since the set of all finite unions of finite intersections of the elements of the family  $\PP\cup\PP^*$ coincides with $\BB$ (see, e.g., \cite{Sik}), we get that  the set of all finite unions of  finite intersections of the elements of the family $\QQ\cup\QQ^*$ coincides with $CO(Y)$. Since $CO(Y)$ is a base of $Y$ (by the Stone Duality Theorem \cite{ST}), we get that $\QQ\cup\QQ^*$ is a subbase of $Y$. Therefore, $\QQ$ is an almost subbase of $Y$. So, $Y$ is a zero-dimensional EC.

The proof for the case when $\PP$ is a $\sigma$-boundedly point-finite family is analogous.
\sqs

The characterization of the subspaces of ECa, UECa and SECa, given by Theorem \ref{MainE}, is in the spirit of the Rosenthal characterization of ECa \cite{Ro} (see Theorem \ref{Rosenthal} here). The Arhangel'ski\v{\i} Question \ref{p2014} shows that it is desirable to obtain some new  characterizations of the subspaces of ECa, UECa and SECa which are in the spirit of the Michael-Rudin characterization of ECa \cite{MR1} (see Theorem \ref{MR1Th14} here). The next observation is in this direction.

 \begin{theorem}\label{th85} A compact Hausdorff space is an EC iff it has an F-separating
$\sigma$-strongly point-finite family of open subsets.
\end{theorem}

\doc ($\Rightarrow$) In the proof of the necessity of \cite[Theorem 3.4, p.76]{D4} it is noted that every EC has an F-separating
$\sigma$-strongly point-finite family of open subsets.

\smallskip

\noindent($\Leftarrow$) This follows from the Michael-Rudin Theorem \ref{MR1Th14}.
\sqs

The next question arises naturally from  Question \ref{p2014} and Theorem \ref{c87}.

\begin{question}\label{p81}
\rm
Is it true that if a Tychonoff space $X$ has a $\s$-disjoint base then it has a $\s$-disjoint base consisting of cozero-sets?
\end{question}

Let's note that in order to answer in the affirmative Question \ref{p2014} one needs to show, according to Theorem \ref{MainE}, that every Tychonoff space with a $\s$-disjoint base  has a $\s$-strongly point-finite almost subbase.

We will give below some results which are in the spirit of the Question \ref{p81}.

\begin{pro}\label{th86} Every open subset of an E-space is the union of a $\sigma$-point-finite
collection of cozero sets.
\end{pro}

\doc As it is noted on p. 493 (after  Theorem 5.1) of the Michael-Rudin paper \cite{MR1}, {\em every open subset of an EC is the union of a $\s$-point-finite collection of cozero-sets}. Hence, applying Theorem \ref{MainE}, we complete the proof of our assertion.
\sqs

\begin{cor}\label{c89} If an E-space $(X,\TT)$ has a
$\sigma$-point-finite base, then $X$ has a $\sigma$-point-finite base consisting of cozero sets.
\end{cor}

\doc Let $\BB=\bigcup\{\BB_n\st n\in\NNNN\}$ be a base of $X$ and, for every $n\in\NNNN$, $\BB_n$ be a point-finite family.

Let $n\in\NNNN$. Then, for every $U\in\BB_n$, we have (by Proposition \ref{th86}) that $U=\bigcup\AA_U$, where $\AA_U=\bigcup\{\AA_U^m\st m\in\NNNN\}$ and, for each $m\in\NNNN$, $\AA_U^m$ is a point-finite family of cozero-sets.

For every $m,n\in\NNNN$, set
$$\AA_n^m=\bigcup\{\AA_U^m\st U\in\BB_n\}.$$
We will show that for every $m,n\in\NNNN$, $\AA_n^m$ is a point-finite family and
$$\AA=\bigcup\{\AA_n^m\st m,n\in\NNNN\}$$
is a base of $X$.

Let $x\in X$ and $V$ be an open neighborhood of $x$. Then there is $U\in\BB$ such that $x\in U\sbe V$. There exists $n\in\NNNN$ such that $U\in\BB_n$. Then, since $U=\bigcup\AA_U$, there exists $W\in\AA_U$ such that $x\in W\sbe U$. Further, there exists $m\in\NNNN$ such that $W\in\AA_U^m$. Then $W\in\AA_n^m\sbe \AA$ and $x\in W\sbe V$. So, $\AA$ is a base of $X$ and $\AA\sbe Coz(X)$.

Let $m,n\in\NNNN$. We will show that $\AA_n^m$ is a point-finite family. Let $x\in X$. If $x\nin\bigcup\BB_n$ then $x\nin\bigcup\AA_n^m$ and all is clear. Suppose that $x\in\bigcup\BB_n$. Since $\BB_n$ is a point-finite collection, there exist only finitely many elements $U_1,\ldots,U_k$ of $\BB_n$ such that $x\in U_i$ for every $i=1,\ldots,k$. Let $i\in\{1,\ldots,k\}$ and $x\in\bigcup\AA^m_{U_i}$. Then $x$ belongs to only finitely many elements of $\AA^m_{U_i}$, because $\AA^m_{U_i}$ is a point-finite family. Therefore, the point $x$ can belong to only finitely many elements of the family $\AA^m_n$. Hence, $\AA$ is a $\s$-point-finite family.
\sqs

\begin{cor}\label{c899} If the answer to the Question \ref{p2014} is "Yes", then every Tychonoff space with a
$\sigma$-disjoint base has a $\sigma$-point-finite base consisting of cozero sets.
\end{cor}

\begin{pro}\label{th88} Every open subset of a UE-space $(Y,\TT)$ is the union of a $\sigma$-boun\-dedly point-finite family of cozero sets.
\end{pro}

\doc It follows from the proof of Theorem \ref{th87} (exactly as Proposition \ref{th86} was derived from the proof of Theorem \ref{MR1Th14}). Indeed, let $V\in\TT$. Using Theorem \ref{MainE}, we get that $Y$ has a compactification $X$ which is a UEC. Let $U=Ex_{X}(V)$. Then, in the notation of the proofs of the sufficiency  of Theorems  \ref{MR1Th14} and \ref{th87}, we can suppose that $U\in\UU$. Further, we have that $\{W_{\a,m}(U)\st \a<\l, m\in\NNNN\}$ is an open cover of $U$ (because $S_{\a,m}(U)\sbe W_{\a,m}(U)\sbe U$ and $\{S_{\a,m}(U)\st \a<\l, m\in\NNNN\}$ covers $U$) (see the proof of the sufficiency of Theorem  \ref{MR1Th14}). 
Then, clearly, $V=\bigcup\{Y\cap W_{\a,m}(U)\st \a<\l, m\in\NNNN\}$.  
Since $\{W_{\a,m}(U)\st \a<\l, m\in\NNNN\}\sbe \WW$ and $\WW$ is a $\sigma$-boundedly point-finite family of cozero sets in $X$ (as it is shown in the proof of the sufficiency of Theorem \ref{th87}),
we get that 
$V$ is a union of a $\sigma$-boundedly point-finite family of cozero subsets of $Y$.
\sqs

\begin{cor}\label{c810}  If $X$ has a $\sigma$-disjoint base and $X$ is a UE-space then it has a
$\sigma$-strongly point-finite base consisting of cozero-sets.
\end{cor}

\doc It is similar to the proof of Corollary \ref{c89}.

Let $\BB=\bigcup\{\BB_n\st n\in\NNNN\}$ be a base of $X$ and, for every $n\in\NNNN$, $\BB_n$ be a disjoint family.

Let $n\in\NNNN$. Then, for every $U\in\BB_n$, we have (by Proposition \ref{th88}) that $U=\bigcup\AA_U$, where $\AA_U=\bigcup\{\AA_U^m\st m\in\NNNN\}$ and, for each $m\in\NNNN$, $\AA_U^m$ is a boundedly point-finite family of cozero-sets.

For every $m,n\in\NNNN$, set
$$\AA_n^m=\bigcup\{\AA_U^m\st U\in\BB_n\}.$$
We will show that for every $m,n\in\NNNN$, $\AA_n^m$ is a strongly point-finite family and
$$\AA=\bigcup\{\AA_n^m\st m,n\in\NNNN\}$$
is a base of $X$.

The fact that $\AA$ is a base can be established exactly as in the proof of Corollary \ref{c89}.

Let $m,n\in\NNNN$. We shall show that $\AA_n^m$ is a strongly point-finite family.
Let $\mu\sbe \AA_n^m$ and $|\mu|=\aleph_0$. Let us set $\mu=\{V_k\st k\in\NNNN\}$. Then, for every $k\in\NNNN$, there exists $U_k\in\BB_n$ such that $V_k\in\AA_{U_k}^m$.
Hence, $V_k\sbe U_k$, for every $k\in\NNNN$. Now, we can define a function $\p:\NNNN\lra\BB_n$, $k\mapsto U_k$. There are two possibilities:

\smallskip

\noindent(a) {\em $\p$ is not a constant function.} Then there exist $k,l\in\NNNN$ such that $U_k\neq U_l$. Since $\BB_n$ is a disjoint family, we get that $U_k\cap U_l=\ems$. Hence $V_k\cap V_l=\ems$ and we set $\mu\ap=\{V_k,V_l\}$. Then $\mu\ap$ is a finite subfamily of $\mu$ and $\bigcap\mu\ap=\ems$.

\smallskip

\noindent(b) {\em $\p$ is a constant function.} Then $U_k=U_1$ for every $k\in\NNNN$. Set $U=U_1$. Then $\mu\sbe \AA_U^m$. Since the order of the family $\AA_U^m$ is finite, we get that the order of the family $\mu$ is finite. Thus there exists a finite subfamily $\mu\ap$ of $\mu$ such that $\bigcap\mu\ap=\ems$.

Therefore, $\AA_n^m$ is a strongly point-finite family. Thus  $\AA$ is a $\s$-strongly point-finite base of $X$ consisting of cozero-sets.
\sqs

\bigskip

\noindent{\bf Acknowledgements.} The author would like to thank Prof. A. V. Arhangel'ski\v{\i} for the useful discussions on the subject.


\baselineskip = 0.75\normalbaselineskip

\end{document}